\documentclass[a4paper,11pt]{article}
\usepackage{fullpage}
\usepackage{amsmath}
\usepackage{amssymb}
\usepackage{siunitx}
\usepackage{setspace}
\usepackage[pdftex]{graphicx}
\usepackage[usenames, dvipsnames]{color}                        %%% add "monochrome" to the Options to avoid colors
\usepackage{booktabs}
\usepackage{multirow}
\usepackage{subfigure}
\usepackage{authblk}
\usepackage{natbib}                                 %%% use \citep[][]{} or \citet[][]{}
\usepackage{palatino}
\usepackage{textcomp} %for \textpm

\usepackage{hyperref}                                           %%% golden rule: should be loaded as the last package
\hypersetup{%                                                   %%% see e.g., http://www.ce.cmu.edu/~kijoo/latex2pdf.pdf ; http://www.tug.org/applications/hyperref/
bookmarksnumbered,
pdfstartview={FitH},
linkcolor=blue,
citecolor=blue,
colorlinks=true,
pdfpagemode={UseNone},
plainpages=false
}%
\newcommand{\cV}{\mathcal{V}}
\newcommand{\cA}{\mathcal{A}}
\newcommand{\cK}{\mathcal{K}}

\newcommand{\cU}{\mathcal{U}}
\newcommand{\cR}{\mathcal{R}}

\allowdisplaybreaks

\title{A Comparison of Models for Uncertain Network Design}
\author[1]{Francis Garuba\thanks{Corresponding author. Email: \texttt{f.garuba@lancaster.ac.uk}}} \author[2]{Marc Goerigk} \author[1]{Peter Jacko}
\affil[1]{{\small Department of Management Science, Lancaster University, United Kingdom}}
\affil[2]{{\small Network and Data Science Management, University of Siegen, Germany}}

\date{}

\begin{document}

\maketitle

\begin{abstract}
To solve a real-world problem, the modeler usually needs to make a trade-off between model complexity and usefulness. This is also true for robust
optimization, where a wide range of models for uncertainty, so-called uncertainty sets, have been proposed. However, while these sets have been mainly
studied from a theoretical perspective, there is little research comparing different sets regarding their usefulness for a real-world problem.

In this paper we consider a network design problem in a telecommunications context. We need to invest into the infrastructure, such that there is
sufficient capacity for future demand which is not known with certainty. There is a penalty for an unsatisfied realized demand, which needs to be
outsourced. We consider three approaches to model demand: using a discrete uncertainty set, using a polyhedral uncertainty set, and using the mean of a
per-commodity fitted zero-inflated uniform distribution. While the first two models are used as part of a robust optimization setting, the last model
represents a simple stochastic optimization setting. We compare these approaches on an efficiency frontier real-world data taken from the online library
SNDlib and observe that, contrary to current research trends, robust optimization using the polyhedral uncertainty set may result in less efficient
solutions.
\end{abstract}


\section{Introduction}

Network design models have found wide application in the planning, design and operations management of transportation, power \& energy distribution, supply chain logistic and telecommunications networks. Usually, they are based on mixed-integer programming models, and many such models have been developed over the decades for network design and expansion problems, see, e.g.,  \cite{Magnanti1984,Minoux1989,Bertsekas1998}.

In telecommunications for instance, network design models can be used to curb congestion and to provide an acceptable quality of service to the subscribers. Effort to provide an acceptable service has resulted in capital expenditure of billions of USD in global telecoms investment. Optimization of investments has thus attained a key strategic role in this industry. Moreover, these decisions need to be made well ahead of time based on a forecast of future traffic demand.

Unfortunately, traffic demand has proven to be difficult to predict accurately. In order to factor in this uncertainty and design a network that is immune to traffic variability, robust optimization approaches have been proposed.
For this purpose, a number of uncertainty models have already been developed and investigated (see \cite{goerigk2016algorithm,Ben-Tal2009,Bertsimas2011}). The drawback of classic approaches, however, is that
the uncertainty set is assumed to be given, i.e., the decision maker can advise us how the uncertainty is shaped. Moreover, an inappropriate choice of uncertainty set may result in models that are too conservative or in some cases computationally intractable. As the decision maker cannot be expected to make this choice in practice,  data-driven and learning approaches have been recently proposed (see \cite{Bertsimas2017,Chassein2018}).

To the best of our knowledge, we follow this approach for the first time for network design problems, by comparing which uncertainty set actually fits real-world data. We compare two robust optimization approaches for a network capacity expansion model with outsourceable demand (see, e.g., \cite{Bertsekas1998,Bektas2009}). In this setting, we need to invest into the network infrastructure now, so that each commodity can be routed to satisfy its uncertain demand later. Demand which cannot be satisfied is outsourced, which is modeled through a linear penalty on its amount.

The two approaches under consideration are (1) a discrete uncertainty set which assumes that all demands are in closed form; and (2) a polyhedral set with wider range of possible scenarios which results in a heuristic mix-integer program to solve the resulting robust problem. These two are compared on real-world data taken from SNDlib and also compared with performance of a third model outside the robust framework, a simple stochastic optimization approach.

The rest of this paper is organized as follows. \autoref{sec:literature} presents a literature review of related research. In \autoref{sec:problem}, we introduce the problem description of robust network capacity expansion with outsourcing and mathematical models for both the discrete and polyhedral uncertainty sets with detailed construction of the robust counterparts. Experimental results and main findings using data from the SNDlib (see \cite{Orlowski2010}) are discussed in \autoref{sec:computational}. Finally, \autoref{sec:conclusion} concludes our work and points out future research directions.

\section{Literature Review }
\label{sec:literature}

The study of uncertainties in decision problems has resulted in two broad areas of research, namely \emph{stochastic} (see, e.g., \cite{Birge2011}) and \emph{robust} (see, e.g., \cite{Ben-Tal2009}) optimization frameworks. While the stochastic approach usually assumes that a probability distribution of the uncertain data is known with precision, the robust approach assumes that the uncertain data lies within a predetermined set. The renewed interest in the latter can be attributed to the works of \cite{Ben_Tal_1999} and \cite{El_Ghaoui_1998} with many other collaborators.

The two frameworks also have a dynamic context, where a part of the decision has to be made after the realization of the uncertain data. This is known as two-stage stochastic and robust optimization. Depending on the context, two-stage robust problems are also known as adjustable robust counterparts (ARC). Here, the decision variables are partitioned into two sets: the non-adjustable ones ("here and now" decisions) which must be fixed in advance before the realization of the uncertainty sets and the adjustable ones ("wait and see" decisions), which are computed after the uncertain parameters are revealed (\cite{Ben_Tal_2004}).

As the ARC is more representative of real life situations where decisions are made over multiple periods, this framework has attracted interest from the research community. However, its its general form is known to be computationally intractable, which has led to an approximate model using affine decision rules \textit{(ADR)}. In this affine adjustable robust counterpart \textit{(AARC)}, the adjustable part of the decision is assumed to be an affine linear function of the uncertain data (\cite{Ben_Tal_2004}). This emulates a linear feedback as a controller to adjust for the desired output.

Just like in many other fields, robust optimization has found increasing use and application in the network design area. \cite{Atamt_rk_2007} considered a two-stage robust network flow problem under demand uncertainty following the work of \cite{Ben_Tal_2004}, while \cite{Ouorou2007} introduced affine routing in the their robust network capacity planning model. \cite{Ord_ez_2007} looked at network capacity expansion under both demand and cost uncertainty. \cite{Koster_2013} considered a robust network design problem with static routing in the setting of \cite{Bertsimas_2004}. \cite{Poss_2012} apply the AARC to robust network design with polyhedral uncertainty and \cite{Babonneau2013} used a refined version of ADR in their robust capacity assignment for networks with uncertain demand.  Recently, \cite{Pessoa_2015} used a cutting plane algorithm while taking into consideration the uncertainty in unmet demand outsourced cost.

Regarding uncertainty sets, polyhedral sets are most frequently used in radio network design, along with hose models from the works of \cite{Duffield1999,Fingerhut_1997}, budget uncertainty by \cite{Atamt_rk_2007} and cardinal constrained uncertainty by \cite{Bertsimas_2004}, and interval uncertainty among others.

Little research compares these models. \cite{Atamt_rk_2007} compared their single-stage robust model using budget uncertainty with a scenario-based
two-stage stochastic approach. \cite{Chassein2018} constructed different uncertainty sets from real world data and compared performance within and outside
sample for shortest path problems. Our focus is to compare the discrete and the polyhedral uncertainty sets in network capacity expansion, to arrive at
which one better fits real-world data, while also comparing to the performance of a simple stochastic model using the mean demand.

\section{Problem Description}
\label{sec:problem}

We consider a multi-commodity network flow design problem where incremental capacities are installed in response to uncertain traffic demand. The problem is modeled in a way that allows for capacity expansion such that routing of traffic for the different commodities over the arcs subject to design and network constraints is possible while minimizing the total cost involved. We refer to this model as the robust network capacity expansion problem \textit{(RNCEP)}.

\subsection{The Basic RNCEP}
\label{sec:basic problem}

The network under consideration can be represented by a directed graph, $G=(\cV, \cA)$. Each of the arcs $a \in \cA$ has an original capacity $u_a$. The original capacity on each arc $ a $ can be upgraded at a cost $c_a$ per each additional unit $x_a$ of capacity. There is a set of commodities $\cK=\{1,\ldots,K\}=:[K]$ which need to be routed across the network, each commodity $k \in \cK$ consisting of a demand $d^k \ge 0$, a source node $s^k\in\cV$, and a sink node $t^k\in\cV$. Additionally, let $\sigma$ be the cost of not satisfying one unit of demand over the planning horizon (i.e., by outsourcing it). If all demands are known, the nominal network capacity expansion problem can then be formulated as follows:
\begin{align}
\min\ &\sum_{a\in \cA} c_a x_a + \sigma \sum_{k\in\cK} \left[ d_k - \sum_{a\in\delta^-(t^k)} f^k_a + \sum_{a\in\delta^+(t^k)} f^k_a \right]_+ \label{con1}\\
\text{s.t. } &  \sum_{a\in \delta^-(v)} f^k_a - \sum_{a\in \delta^+(v)} f^k_a \ge
0
& \forall k\in \cK,  v\in \cV\setminus\{s^k,t^k\} \label{con2} \\
& \sum_{k\in \cK} f^k_a \leq u_a + x_a & \forall   a\in \cA  \label{con3}\\
& f^k_a \ge 0 & \forall k\in\cK,d\in\cU,a\in\cA \label{con4} \\
& x_a \geq 0 & \forall a\in\cA \label{con5}
\end{align}
Here, $[y]_+$ denotes $\max\{0,y\}$, while $\delta^+(v)$ and $\delta^-(v)$ are the sets of the outgoing and incoming arc at node $v \in \cV$, respectively. Variables $f^k_a$ denote the flow of commodity $k\in\cK$ along edge $a\in\cA$, while $x_a$ models the amount of capacity being added to arc $a$. The objective function~\eqref{con1} is to minimize the sum of capacity expansion cost and outsourcing costs. Constraints~\eqref{con2} are a variant of flow constraints, where we allow an arbitrary amount of flow to leave the source node $s^k$. Through the objective, only the flow arriving in $t^k$ is counted. It is allowed to diminish the flow outside of $s^k$ and $t^k$; note that there is an optimal solution where this does not happen. We do not assume equality in Constraints~\eqref{con2} to apply our robust optimization approach in the following section. Finally, Constraints~\eqref{con3} model the capacity on each edge.

The actual demand values $\pmb{d}$ are uncertain, and can take any value in a predetermined uncertainty set $\cU$. The two sets under consideration in this work are the \textit{discrete uncertainty set}, which can be represented as  $\cU = \{\pmb{d}^1,\ldots,\pmb{d}^N\}$, and the \textit{polyhedral uncertainty set}, which can be represented as $\cU = \left\{ \pmb{d}\in\mathbb{R}^K_+ : V\pmb{d} \le \pmb{b}, d_k\in[\underline{d}_k,\overline{d}_k] \right\}$.

The robust network capacity expansion problem then is to find a minimum installation cost of additional capacities while satisfying all potential traffic demands such that actual flows do not exceed cumulative link capacities whatever the realization of demands in $\cU$. Thus, the RNCEP is a two stage robust problem with recourse applying the general framework of \cite{Ben_Tal_2004}. The capacity expansion represented by variables $\pmb{x}$ is the first stage decision variable which has to be fixed before the realization of $\pmb{d} \in \cU$. Once the uncertain  demand data is revealed, the traffic adjustment takes place by routing a multi-commodity flow with second stage variable $f^k_a(\pmb{d})$. This can be modeled as follows:
\begin{align}
\min\ &\sum_{a\in \cA} c_a x_a + \max_{d\in\cU} \sigma \sum_{k\in\cK} \left[ d_k - \sum_{a\in\delta^-(t^k)} f^k_a(\pmb{d}) + \sum_{a\in\delta^+(t^k)} f^k_a(\pmb{d}) \right]_+ \label{con1a}\\
\text{s.t. } &  \sum_{a\in \delta^-(v)} f^k_a(\pmb{d}) - \sum_{a\in \delta^+(v)} f^k_a(\pmb{d}) \ge
0
& \forall k\in \cK, \pmb{d} \in \cU, v\in \cV\setminus\{s^k,t^k\} \label{con2a} \\
& \sum_{k\in \cK} f^k_a(\pmb{d}) \leq u_a + x_a & \forall  \pmb{d} \in \cU, a\in \cA  \label{con3a}\\
& f^k_a(\pmb{d}) \ge 0 & \forall k\in\cK, \pmb{d}\in\cU,a\in\cA \label{con4a} \\
& x_a \geq 0 & \forall a\in\cA \label{con5a}
\end{align}
Here, we have modified Constraints~(\ref{con1}-\ref{con5}) to take all scenarios into account. Being a robust model, we consider the worst-case costs in Objective~\eqref{con1a}, while all constraints need to hold for all scenarios $\pmb{d}\in\cU$. In the following, we reformulate the general model~(\ref{con1a}-\ref{con5a}) for specific uncertainty sets.

\subsection{Robust Optimization with Discrete Uncertainty}
\subsubsection{Model}

Let $\cU=\{\pmb{d}^1, \dots, \pmb{d}^N\}$ be a discrete uncertainty set, where $N$ is the number of scenarios.
In this case, variables $f^k_a(\pmb{d})$ become $f^{k,i}_a$ for all $i\in[N]$. The robust objective function~\eqref{con1a} is reformulated using additional variables $h^{k,i}:=[d^i_k - \sum_{a\in\delta^-(t^k)} f^{k,i}_a + \sum_{a\in\delta^+(t^k)} f^{k,i}_a]_+$  for $k\in\cK$, $i\in[N]$, and $\tau := \max_{i\in[N]} \sum_{k\in\cK} h^{k,i}$.
The problem then becomes:
\begin{align}
\min\ &\sum_{a\in \cA} c_a x_a + \sigma \tau \label{disct1}\\
\text{s.t. } & \tau \ge \sum_{k\in\cK} h^{k,i} & \forall i\in[N] \label{disct3}\\
& h^{k,i} \ge  d^i_k - \sum_{a\in\delta^-(t^k)} f^{k,i}_a + \sum_{a\in\delta^+(t^k)} f^{k,i}_a & \forall i\in[N], k\in\cK \label{disct4}\\
& \sum_{a\in \delta^-(v)} f^{k,i}_a - \sum_{a\in \delta^+(v)} f^{k,i}_a \ge
0 & \forall k\in \cK, i\in[N], v\in \cV\setminus\{s^k,t^k\} \label{disct2} \\
& \sum_{k\in \cK} f^{k,i}_a \leq u_a + x_a & \forall  i\in[N], a\in \cA \label{disct5} \\
& f^{k,i}_a \ge 0 & \forall k\in\cK,i\in[N],a\in\cA \label{disct6}\\
& h^{k,i} \ge 0 & \forall k\in\cK,i\in[N] \label{disct7}\\
& x_a \geq 0 & \forall a\in\cA
\end{align}
Here, Constraints~\eqref{disct2} and~\eqref{disct5} correspond to Constraints~\eqref{con2a} and~\eqref{con3a}, whereas the additional Constraints~\eqref{disct3} and~\eqref{disct4} are used to ensure variables $\tau$ and $h^{k,i}$ have the intended effect. Note that, as we minimize, the maximum operator can be expressed by using $\ge$-constraints over the set.

\subsubsection{Constructing Data-Based Discrete Uncertainty}
To construct discrete uncertainties uncertainties, we assume that scenarios
\[ \cR = \{ \pmb{r}^1, \ldots, \pmb{r}^N\} \]
of real demands with $\pmb{r}^i \in\mathbb{R}^K_+$ are given, along with the respective source and sink nodes. The trivial approach would be to use directly $\cU=\cR$. However, previous research (see \cite{Chassein2018}) has shown that this may result in an overfitting to the available data. Instead, we consider different scalings. For a fixed commodity $k\in\cK$, let $N'\le N$ denote the absolute frequency that $r^{i,k} > 0$ over all $i\in[N]$. Then
\[ \hat{r}^k = \frac{1}{N'}\sum_{i\in[N]} r^{i,k} \]
be the average of the demand scenarios for each $k\in\cK$. For a given $\lambda\in[0,1]$, we set
$d^{i,k}(\lambda) = \lambda r^{i,k} + (1-\lambda) \hat{r}^k$ and
\[\cU(\lambda) = \left\{ \pmb{d}^1(\lambda), \ldots, \pmb{d}^N(\lambda) \right\}. \]
The case $\lambda = 0$ means that we ignore uncertainty and use the average case, while $\lambda=1$ uses the original demand scenarios $\cR$.

\subsection{Robust Optimization with Polyhedral Uncertainty}
\label{sec:polyhedral}
\subsubsection{Model}

We now assume the demand uncertainty is given through a general polyhedron of the form
\[ \cU = \big\{ \pmb{d}\in\mathbb{R}^K_+ : V\pmb{d} \le \pmb{b}, d_k\in[\underline{d}_k,\overline{d}_k] \big\} \]
where $V=(v_{ik})$ is a matrix in $\mathbb{R}^{M\times K}$ and $\pmb{b}$ is a vector in $\mathbb{R}^{M}$ (i.e., there are $M$ linear constraints on the demand vector). To find a tractable robust counterpart, we apply the framework of affine decision rules (ADR) by restricting the flow variables to be affine functions of the uncertainty, i.e.,
\[ f^k_a(\pmb{d}) = \phi^k_a + \sum_{\ell\in\cK} \Phi^{k,\ell}_a d_{\ell} \]
with $\phi^k_a$ and $\Phi^{k,\ell}_a$ being unknown coefficients of the affine linear function in $\pmb{d}$. We now consider each constraint and the objective of problem~(\ref{con1a}-\ref{con5a}) and reformulate them using strong duality.

By substituting for $f^k_a(\pmb{d})$, the flow constraints~\eqref{con2a} become:
\[ \sum_{a\in\delta^-(v)} \left( \phi^k_a + \sum_{\ell\in\cK} \Phi^{k,\ell}_a d_{\ell} \right) - \sum_{a\in\delta^+(v)} \left(\phi^k_a + \sum_{\ell\in\cK} \Phi^{k,\ell}_a d_{\ell}\right) \ge 0 \qquad \forall k\in\cK, v\in\cV\setminus\{s^k,t^k\}, \pmb{d}\in\cU, \]
which is equivalent to
\begin{equation}
\sum_{a\in\delta^-(v)}\phi^k_a - \sum_{a\in\delta^+(v)} \phi^k_a \ge \sum_{\ell\in\cK} \left(\sum_{a\in\delta^+(v)}\Phi^{k,\ell}_a - \sum_{a\in\delta^-(v)} \Phi^{k,\ell}_a\right) d_\ell \qquad \forall k\in\cK, v\in\cV\setminus\{s^k,t^k\}, \pmb{d}\in\cU.\label{con2b}
\end{equation}
For each $k\in\cK, v\in\cV\setminus\{s^k,t^k\}$ we can write the worst-case problem as
\begin{align*}
\max\  & \sum_{\ell\in\cK} \left(\sum_{a\in\delta^+(v)}\Phi^{k,\ell}_a - \sum_{a\in\delta^-(v)} \Phi^{k,\ell}_a\right) d_\ell \\
\text{s.t. } & \sum_{\ell\in\cK} v_{i\ell} d_\ell \le b_i & \forall i\in[M] && [\alpha^{k,v}_i] \\
& d_{\ell} \le \overline{d}_\ell & \forall \ell\in\cK && [ \overline{\beta}^{k,v}_\ell]\\
& -d_{\ell} \le -\underline{d}_{\ell} & \forall \ell\in\cK && [\underline{\beta}^{k,v}_{\ell}]
\end{align*}
We now consider the dual of this linear optimization problem. In brackets behind every constraint of the primal problem, we have listed the corresponding dual variable. The dual problem then becomes
\begin{align*}
\min\ & \sum_{i\in[M]} b_i\alpha^{k,v}_i + \sum_{\ell\in\cK} ( \overline{d}_l\overline{\beta}^{k,v}_{\ell} - \underline{d}_\ell\underline{\beta}^{k,v}_\ell) \\
\text{s.t. } &  \sum_{i\in[M]} v_{i\ell}\alpha^{k,v}_i + \overline{\beta}^{k,v}_\ell - \underline{\beta}^{k,v}_\ell \ge \sum_{a\in\delta^+(v)}\Phi^{k,\ell}_a - \sum_{a\in\delta^-(v)} \Phi^{k,\ell}_a  &\forall \ell\in\cK \\
& \alpha^{k,v}_i \ge 0 & \forall i\in[M] \\
&  \overline{\beta}^{k,v}_\ell \ge 0 & \forall \ell\in\cK \\
& \underline{\beta}^{k,v}_{\ell} \ge 0 & \forall \ell \in\cK.
\end{align*}
By applying strong duality, we can conclude that the optimal objective value of this dual problem is equal to the worst-case of the right-hand side of Constraint~\eqref{con2b}.

Overall, \textbf{Constraint~\eqref{con2a}} is replaced by the following set of constraints and variables:
\begin{align*}
& \sum_{a\in\delta^-(v)}\phi^k_a - \sum_{a\in\delta^+(v)} \phi^k_a \ge \sum_{i\in[M]} b_i\alpha^{k,v}_i + \sum_{\ell\in\cK} ( \overline{d}_l\overline{\beta}^{k,v}_{\ell} - \underline{d}_\ell\underline{\beta}^{k,v}_\ell) & \forall  k\in\cK, v\in\cV\setminus\{s^k,t^k\} \\
& \sum_{i\in[M]} v_{i\ell}\alpha^{k,v}_i + \overline{\beta}^{k,v}_\ell - \underline{\beta}^{k,v}_\ell \ge \sum_{a\in\delta^+(v)}\Phi^{k,\ell}_a - \sum_{a\in\delta^-(v)} \Phi^{k,\ell}_a  &\forall  k,\ell\in\cK, v\in\cV\setminus\{s^k,t^k\} \\
& \alpha^{k,v}_i \ge 0 & \forall i\in[M], k\in\cK, v\in\cV\setminus\{s^k,t^k\} \\
&  \overline{\beta}^{k,v}_\ell \ge 0 & \forall k,\ell\in\cK, v\in\cV\setminus\{s^k,t^k\}\\
& \underline{\beta}^{k,v}_{\ell} \ge 0 & \forall k,\ell \in\cK, v\in\cV\setminus\{s^k,t^k\}
\end{align*}
We follow a similar procedure for the other constraints. Constraint~\eqref{con3a} can be rewritten as
\[ \sum_{k\in\cK} \left( \phi^k_a + \sum_{\ell\in\cK} \Phi^{k,\ell}_a d_{\ell} \right) \le u_a + x_a \qquad \forall d\in\cU, a\in\cA \]
The subproblem
\begin{align*}
\max\ &\sum_{\ell\in\cK} (\sum_{k\in\cK} \Phi^{k,\ell}_a ) d_\ell \\
\text{s.t. } & \pmb{d}\in\cU
\end{align*}
has the same structure as before. Using dual variables $\pi^a_i,\overline{\rho}^a_\ell,\underline{\rho}^a_\ell$, we can replace \textbf{Constraint~\eqref{con3a}} with the following:
\begin{align*}
& \sum_{k\in\cK} \phi^k_a + \sum_{i\in[M]} b_i \pi^a_i + \sum_{\ell\in\cK} (\overline{d}_\ell\overline{\rho}^a_\ell - \underline{d}_{\ell}\underline{\rho}^a_\ell) \le u_a + x_a & \forall a\in\cA \\
& \sum_{i\in[M]} v_{i\ell}\pi^a_i + \overline{\rho}^a_\ell -\underline{\rho}^a_\ell \ge \sum_{k\in\cK} \Phi^{k,\ell}_a & \forall \ell\in\cK, a\in\cA \\
& \pi^a_i \ge 0 & \forall i\in[M],a\in\cA \\
& \overline{\rho}^a_\ell \ge 0 & \forall a\in\cA,\ell\in\cK \\
& \underline{\rho}^a_\ell \ge 0 & \forall a\in\cA,\ell\in\cK
\end{align*}
We now consider the positivity constraint~\eqref{con4a}. This becomes
\[ \phi^k_a + \sum_{\ell\in\cK} \Phi^{k,\ell}_a d_{\ell} \ge 0 \qquad \forall k\in\cK,a\in\cA,\pmb{d}\in\cU \]
Using duality with variables $\xi^{k,a}_i,\overline{\zeta}^{k,a}_\ell,\underline{\zeta}^{k,a}_\ell$ we replace \textbf{Constraint~\eqref{con4a}} with the following:
\begin{align*}
& \phi^k_a \ge \sum_{i\in[M]} b_i\xi^{k,a}_i + \sum_{\ell\in\cK} (\overline{d}_\ell \overline{\zeta}^{k,a}_\ell - \underline{d}_\ell \underline{\zeta}^{k,a}_\ell) & \forall k\in\cK,a\in\cA \\
& \sum_{i\in[M]} v_{i\ell} \xi^{k,a}_i + \overline{\zeta}^{k,a}_\ell - \underline{\zeta}^{k,a}_\ell \ge - \Phi^{k,\ell}_a & \forall k,\ell\in\cK,a\in\cA \\
& \xi^{k,a}_i \ge 0 & \forall k\in\cK,a\in\cA,i\in[M] \\
&\overline{\zeta}^{k,a}_\ell \ge 0 & \forall k,\ell\in\cK,a\in\cA \\
&\underline{\zeta}^{k,a}_\ell \ge 0 & \forall k,\ell\in\cK,a\in\cA
\end{align*}
Finally, we consider the objective function~\eqref{con1a}. We need to solve the following problem:
\begin{align*}
\max\ &\sum_{k\in\cK} \left[ d_k - \sum_{a\in\delta^-(t^k)} \left(\phi^k_a + \sum_{\ell\in\cK} \Phi^{k,\ell}_a d_\ell \right) +  \sum_{a\in\delta^+(t^k)} \left(\phi^k_a + \sum_{\ell\in\cK} \Phi^{k,\ell}_a d_\ell \right) \right]_+ \\
\text{s.t. } & \sum_{\ell\in\cK} v_{i\ell} d_\ell \le b_i & \forall i\in[M] \\
& d_\ell \le \overline{d}_\ell & \forall \ell\in\cK \\
& -d_\ell \le -\underline{d}_\ell & \forall \ell \in\cK
\end{align*}
We introduce new variables $z_k\in\{0,1\}$ to remove the positivity bracket from the objective.
\begin{align*}
\max \ & \sum_{k\in\cK} \left( d_k - \sum_{a\in\delta^-(t^k)} \left(\phi^k_a + \sum_{\ell\in\cK} \Phi^{k,\ell}_a d_\ell \right) +  \sum_{a\in\delta^+(t^k)} \left(\phi^k_a + \sum_{\ell\in\cK} \Phi^{k,\ell}_a d_\ell \right) \right) z_k \\
\text{s.t. } & \sum_{\ell\in\cK} v_{i\ell} d_\ell \le b_i & \forall i\in[M] \\
& d_\ell \le \overline{d}_\ell & \forall \ell\in\cK \\
& -d_\ell \le -\underline{d}_\ell & \forall \ell \in\cK\\
& z_k\in\{0,1\} & \forall k\in\cK
\end{align*}
We set $z'_{k,\ell} := d_\ell z_k$ and get
\begin{align*}
\max\ & \sum_{k\in\cK} \left( z'_{kk} - \sum_{a\in\delta^-(t^k)} \left(\phi^k_az_k + \sum_{\ell\in\cK} \Phi^{k,\ell}_a z'_{k\ell} \right) +  \sum_{a\in\delta^+(t^k)} \left(\phi^k_az_k + \sum_{\ell\in\cK} \Phi^{k,\ell}_a z_{k\ell} \right) \right) \\
\text{s.t. } & \sum_{\ell\in\cK} v_{i\ell} d_\ell \le b_i & \forall i\in[M] && [\mathfrak{q}_i]\\
& z'_{k\ell} \le d_\ell & \forall k,\ell\in\cK && [\mathfrak{r}_{k\ell}]\\
& z'_{k\ell} \le \overline{d}_\ell z_k & \forall k,\ell\in\cK && [\mathfrak{s}_{k\ell}]\\
& d_\ell + \overline{d}_\ell z_k - z'_{k\ell} \le \overline{d}_\ell & \forall k,\ell\in\cK && [\mathfrak{t}_{k\ell}] \\
& d_\ell \le \overline{d}_\ell & \forall \ell\in\cK && [\mathfrak{u}_{\ell}] \\
& -d_\ell \le -\underline{d}_\ell & \forall \ell \in\cK && [\mathfrak{v}_\ell] \\
& z_k\in\{0,1\} & \forall k\in\cK && [\mathfrak{w}_k] \\
& z'_{k\ell} \ge 0 & \forall k,\ell\in\cK
\end{align*}
By relaxing constraints $z_k\in\{0,1\}$ to $z_k\in[0,1]$ for a conservative approximation and dualizing the problem, we arrive at
\begin{align*}
\min\ & \sum_{i\in[M]} b_i\mathfrak{q}_i + \sum_{k\in\cK}\sum_{\ell\in\cK}\overline{d}_\ell\mathfrak{t}_{k\ell} + \sum_{\ell\in\cK} \overline{d}_\ell\mathfrak{u}_\ell - \sum_{\ell\in\cK} \underline{d}_\ell\mathfrak{v}_\ell + \sum_{k\in\cK} \mathfrak{w}_k \\
\text{s.t. } & \sum_{i\in[M]} v_{i\ell} \mathfrak{q}_i - \sum_{k\in\cK} \mathfrak{r}_{k\ell} + \sum_{k\in\cK} \mathfrak{t}_{k\ell} + \mathfrak{u}_\ell - \mathfrak{v}_\ell \ge 0 & \forall \ell\in\cK \\
& -\sum_{\ell\in\cK} \overline{d}_\ell \mathfrak{s}_{k\ell} + \sum_{\ell\in\cK}\overline{d}_\ell \mathfrak{t}_{k\ell} + \mathfrak{w}_k \ge \sum_{a\in\delta^+(t^k)} \phi^k_a - \sum_{a\in\delta^-(t^k)} \phi^k_a & \forall k\in\cK \\
& \mathfrak{r}_{k\ell}  + \mathfrak{s}_{k\ell} - \mathfrak{t}_{k\ell} \ge 1_{k=\ell} + \sum_{a\in\delta^+(t^k)} \Phi^{k,\ell}_a - \sum_{a\in\delta^-(t^k)} \Phi^{k,\ell}_a & \forall k,\ell\in\cK \\
& \mathfrak{q}_i \ge 0 & \forall i\in[M] \\
& \mathfrak{r}_{k\ell},\mathfrak{s}_{k\ell},\mathfrak{t}_{k\ell} \ge 0 & \forall k,\ell\in\cK  \\
& \mathfrak{u}_\ell,\mathfrak{v}_\ell,\mathfrak{w}_\ell \ge 0 & \forall \ell\in\cK
\end{align*}
Overall, we get the following affine adjustable robust counterpart to Problem~(\ref{con1a}-\ref{con5a}):
\begin{align*}
\min\ & \sum_{a\in\cA} c_ax_a + \sigma\left(\sum_{i\in[M]} b_i\mathfrak{q}_i + \sum_{k\in\cK}\sum_{\ell\in\cK}\overline{d}_\ell\mathfrak{t}_{k\ell} + \sum_{\ell\in\cK} \overline{d}_\ell\mathfrak{u}_\ell - \sum_{\ell\in\cK} \underline{d}_\ell\mathfrak{v}_\ell + \sum_{k\in\cK} \mathfrak{w}_k\right)\hspace*{-2.5cm} \\
\text{s.t. } & \sum_{i\in[M]} v_{i\ell} \mathfrak{q}_i - \sum_{k\in\cK} \mathfrak{r}_{k\ell} + \sum_{k\in\cK} \mathfrak{t}_{k\ell} + \mathfrak{u}_\ell - \mathfrak{v}_\ell \ge 0 & \forall \ell\in\cK \\
& -\sum_{\ell\in\cK} \overline{d}_\ell \mathfrak{s}_{k\ell} + \sum_{\ell\in\cK}\overline{d}_\ell \mathfrak{t}_{k\ell} + \mathfrak{w}_k \ge \sum_{a\in\delta^+(t^k)} \phi^k_a - \sum_{a\in\delta^-(t^k)} \phi^k_a & \forall k\in\cK \\
& \mathfrak{r}_{k\ell}  + \mathfrak{s}_{k\ell} - \mathfrak{t}_{k\ell} \ge 1_{k=\ell} + \sum_{a\in\delta^+(t^k)} \Phi^{k,\ell}_a - \sum_{a\in\delta^-(t^k)} \Phi^{k,\ell}_a & \forall k,\ell\in\cK \\
& \sum_{a\in\delta^-(v)}\phi^k_a - \sum_{a\in\delta^+(v)} \phi^k_a \ge \sum_{i\in[M]} b_i\alpha^{k,v}_i + \sum_{\ell\in\cK} ( \overline{d}_l\overline{\beta}^{k,v}_{\ell} - \underline{d}_\ell\underline{\beta}^{k,v}_\ell) & \forall  k\in\cK, v\in\cV\setminus\{s^k,t^k\} \\
& \sum_{i\in[M]} v_{i\ell}\alpha^{k,v}_i +\overline{\beta}^{k,v}_\ell - \underline{\beta}^{k,v}_\ell \ge \sum_{a\in\delta^+(v)}\Phi^{k,\ell}_a - \sum_{a\in\delta^-(v)} \Phi^{k,\ell}_a  &\forall  k,\ell\in\cK, v\in\cV\setminus\{s^k,t^k\} \\
& \sum_{k\in\cK} \phi^k_a + \sum_{i\in[M]} b_i \pi^a_i + \sum_{\ell\in\cK} (\overline{d}_\ell\overline{\rho}^a_\ell - \underline{d}_{\ell}\underline{\rho}^a_\ell) \le u_a + x_a & \forall a\in\cA \\
& \sum_{i\in[M]} v_{i\ell}\pi^a_i + \overline{\rho}^a_\ell - \underline{\rho}^a_\ell \ge \sum_{k\in\cK} \Phi^{k,\ell}_a & \forall \ell\in\cK, a\in\cA \\
& \phi^k_a \ge \sum_{i\in[M]} b_i\xi^{k,a}_i + \sum_{\ell\in\cK} (\overline{d}_\ell \overline{\zeta}^{k,a}_\ell - \underline{d}_\ell \underline{\zeta}^{k,a}_\ell) & \forall k\in\cK,a\in\cA \\
& \sum_{i\in[M]} v_{i\ell} \xi^{k,a}_i + \overline{\zeta}^{k,a}_\ell - \underline{\zeta}^{k,a}_\ell \ge - \Phi^{k,\ell}_a & \forall k,\ell\in\cK,a\in\cA \\
& x_a \ge  0 &\forall a\in\cA \\
& \mathfrak{q}_i \ge 0 & \forall i\in[M] \\
& \mathfrak{r}_{k\ell},\mathfrak{s}_{k\ell},\mathfrak{t}_{k\ell} \ge 0 & \forall k,\ell\in\cK \\
& \mathfrak{u}_\ell,\mathfrak{v}_\ell,\mathfrak{w}_\ell \ge 0 & \forall \ell\in\cK \\
& \alpha^{k,v}_i \ge 0 &  \forall i\in[M], k\in\cK, v\in\cV\setminus\{s^k,t^k\} \\
&  \overline{\beta}^{k,v}_\ell,\underline{\beta}^{k,v}_{\ell} \ge 0 & \forall k,\ell\in\cK, v\in\cV\setminus\{s^k,t^k\}\\
& \pi^a_i \ge 0 & \forall i\in[M],a\in\cA \\
& \overline{\rho}^a_\ell,\underline{\rho}^a_\ell \ge 0 & \forall a\in\cA,\ell\in\cK \\
& \xi^{k,a}_i \ge 0 & \forall k\in\cK,a\in\cA,i\in[M] \\
&\overline{\zeta}^{k,a}_\ell,\underline{\zeta}^{k,a}_\ell \ge 0 & \forall k,\ell\in\cK,a\in\cA
\end{align*}

\subsubsection{Constructing Data-Based Polyhedral Uncertainty}
\label{polcon}

Constructing a polyhedron that contains the demand scenarios $\cR$ can be considered as an optimization problem on its own. We would like to determine constraint coefficients $(v_{i1},\ldots,v_{iK},b_i)$ that determine a polyhedron $\cU$ such that the distance of $\mathcal{R}$ to the boundary of $\cU$ with respect to a norm $\|\cdot\|$ is as small as possible.

Recall that the distance between a point $\pmb{p}$ and a hyperplane $(a_1,\ldots,a_K,b)$ is given through
\[ \frac{ |\sum_{i\in[K]} a_i p_i - b| }{\|\pmb{a}\|^* } \]
where $\|\cdot\|^*$ is the dual norm of $\|\cdot\|$. 
 An optimization model to determine $\cU$ is hence:
\begin{align*}
\min\ & \sum_{i\in[N]} \min_{j\in[M]} (b_j - \sum_{k\in[K]} r^{i,k} v_{jk}) \\
\text{s.t. } & \sum_{k\in[K]} r^{i,k} v_{jk} \le b_j & \forall i\in[N],j\in[M] \\
& \| \pmb{v}_{j\cdot} \|^* = 1
\end{align*}
where $\pmb{v}_{j\cdot}$ denotes the $j$th row of $V$. While such an approach is useful for low-dimensional data (i.e., few commodities $K$), it is less efficient for high-dimensional data. In fact, the additional lower and upper bounds $[\underline{d}_k,\overline{d}_k]$ may already suffice to determine a polyhedron where every point in $\cR$ is on its boundary. Therefore, we also consider randomly generated hyperplanes. To this end, we sample each $v_{ik}$ randomly uniformly from $[0,1]$. Then we set
\[ b_i := \max_{j\in[N]} \sum_{k\in[K]} v_{ik} r^{j,k} \]
to find a tight constraint. In particular, we always contain the sum-constraint where $v_{ik}=1/K$ for all $k\in[K]$.

\subsection{Stochastic Optimization with Distribution Mean}

\subsubsection{Model}
Let $\overline{\pmb{d}}$ be the vector of mean demands of distributions fitted independently to every commodity using demand scenarios $\cR=\{\pmb{r}^1,
\ldots, \pmb{r}^N\}$. We reformulate Problem~(\ref{con1}-\ref{con5}) using only this single mean demand scenario. To linearize the positivity brackets
$[\cdot]_+$, we introduce variables $h^k$ for every commodity $k\in\cK$. The problem then becomes:
\begin{align*}
\min\ &\sum_{a\in \cA} c_a x_a + \sigma \sum_{k\in\cK} h^{k} \\
\text{s.t. } & h^{k} \ge  \overline{d}^k - \sum_{a\in\delta^-(t^k)} f^{k}_a + \sum_{a\in\delta^+(t^k)} f^{k}_a & \forall k\in\cK \\
&  \sum_{a\in \delta^-(v)} f^{k}_a - \sum_{a\in \delta^+(v)} f^{k}_a \ge
0
& \forall k\in \cK, v\in \cV\setminus\{s^k,t^k\}  \\
& \sum_{k\in \cK} f^{k}_a \leq u_a + x_a & \forall a\in \cA \\
& f^{k}_a \ge 0 & \forall k\in\cK,a\in\cA \\
& h^{k} \ge 0 & \forall k\in\cK\\
& x_a \geq 0 & \forall a\in\cA
\end{align*}

\subsubsection{Generating Data-Based Distribution Mean}
\label{sec:stochastic}

The demand for the stochastic optimization model is generated from the demand scenarios $\cR$ using the mean of the zero-inflated uniform distribution in
the following way. For a fixed commodity $k\in\cK$, let $N'\le N$ denote the absolute frequency that $r^{i,k} > 0$ over all $i\in[N]$. To fit a uniform
distribution, set
\[ r^k_{min} = \min_{i\in [N] : r^{i,k} > 0} r^{i,k} \text{ and } r^k_{max}= \max_{i\in [N]}\ r^{i,k}  \]
and the mean of the uniform distribution is $\bar{r}^k=1/2(r^k_{min}+r^k_{max})$. The remaining absolute frequency, $N - N'$, is considered for observing a
zero demand, yielding the mean demand of the zero-inflated uniform distribution $\overline{d}^k = \overline{r}^k N'/N$.

\section{Computational Experiments}
\label{sec:computational}

\subsection{Setup}

The aim of our experiments is to determine which model gives the best solution to uncertain network design. On the one hand, the discrete uncertainty model is simpler than the polyhedral uncertainty model, and we can expect it to be solvable using more commodities, thus giving a more detailed description of the uncertainty. The polyhedral model on the other hand will use less commodities, but has a more complex description of the uncertainty available. As noted in the literature review (Section~\ref{sec:literature}), polyhedral models are popular in current research.

We consider the following experimental setup to address our question. Using a data set of real-world scenarios, we separate it into a training set and an evaluation set. We construct different uncertainty sets only based on the training set, and solve the resulting robust (or stochastic) optimization problems. We then only keep the here-and-now part of the solution, i.e., the decision $\pmb{x}$ on the infrastructure investment. This investment is then assessed on the evaluation set by calculating optimal flows for each scenario. As the first-stage investment costs are already fixed, the flow problem only aims at minimizing the outsourced demand. We then compare investment costs and outsourced demand for all models.

The experimental setup is summarized in \autoref{TabS}. The $\text{Discrete}_1$ experiment fixes two $\sigma$ values for varying values of $\lambda$, while the $\text{Discrete}_2$ experiment fixes $\lambda$ for varying values of $\sigma$. The $\text{Polyhedral}_1$ experiment uses a polyhedron with only one constraint (the sum-constraint) for all eleven values of $\sigma$, $\text{Polyhedral}_2$ uses a polyhedron with two hyperplanes and eight values for $\sigma$, while $\text{Polyhedral}_3$ uses a polyhedron with eight hyperplanes and seven possible $\sigma$ values. The reduced choice for $\sigma$ values with increasing number of hyperplanes was due to increased computation times.

\begin{table}[tbp]\footnotesize
	\centering
	%\tbl{Impact of Choice of k on the presorted original data.}
	\caption{Experimental Setup}
	\begin{tabular}{rcrrrc}
		\toprule
		\textbf{Experiment} & \textbf{Nr of $\sigma$ values} &\textbf{$\sigma$ values} &\textbf{nr $\lambda$ values}&\textbf{$\lambda$ values} & \textbf{Nr of hyperplanes}\\
		\toprule %\hline \hline
		$\text{Discrete}_1$& 2 &    12,450 and 24,900&11&0.0 to 1.0& \\
		$\text{Discrete}_2$& 11 &   0 to 24,900& 2& 0.5 and 1.0&  \\%$ %\widehat{d}$ \& $0.6 \widehat{d}$\\ %\hline
		Stochastic& 11&  0 to 24,900& 1&  1.0& \\
		$\text{Polyhedral}_1$& 11& 0 to 24,900&  &&1\\
		$\text{Polyhedral}_2$& 8 & 0 to 17,430&  &&2\\
		$\text{Polyhedral}_3$& 7 & 0 to 14,940&  &&7\\
		\bottomrule
	\end{tabular}
	%\caption{Experimental setup for generating 120 problem instances for each network.}
	\label{TabS}
\end{table}

The $\text{Discrete}_1$ experiment therefore has to solve $22$ optimization models, and each of these $22$ results was then evaluated on each of the demand scenarios from the evaluation set. The same was carried out for the other five experiments. The choice of $\sigma$, which represent the penalty for unmet demand, was a key consideration for these models and hence in the experimental setup. If $\sigma$ is too small there is incentive for unmet demand where almost all demand are outsourced with no addition of new installed capacity to the network while with a large $\sigma$, the incentive is for negative violation of the constraint which encourages the deployment of new network capacity.

Several values of $\sigma$ were tested in a preliminary experiment using discrete uncertainty, see Table~\ref{Tab2}. Based on the outcomes, the value range for $\sigma$ was selected, taking the $95^{th}$ percentile of the capacity cost distribution into account.

\begin{table}[tbp]\footnotesize
	\centering
	%\tbl{Impact of Choice of k on the presorted original data.}
	\caption{Impact of $\sigma$ on outsourced demand.}
	\begin{tabular}{rcrrrr}
		\toprule
		\textbf{Objective} & \textbf{Commodity} &\textbf{Capacity Add} &\textbf{Sol Time}&\textbf{Outsourced D} & \textbf{Penalty($\sigma$)}\\
		\toprule
		   530,226.88& 400 &       0.00&199.56&127,254.45&   100\\
	 	 5,302,268.77& 400 &       0.00&127.48&127,254.45& 1,000\\%$ %\widehat{d}$ \& $0.6 \widehat{d}$\\ %\hline
		49,191,424.59& 400 &   2,191.83&443.11& 99,409.83&10,000\\
		68,676,799.21& 400 &  16,139.27&372.44& 13,262.58&20,000\\
		72,391,760.98& 400 &  18,151.57&575.90&  6,074.70&30,000\\
		74,113,211.23& 400 &  19,736.41&455.89&  2,335.66&40,000\\
		74,479,094.31& 400 &  20,873.46&378.29&      0.00&50,000\\
		\bottomrule
	\end{tabular}
	%\caption{Experimental setup for generating 120 problem instances for each network.}
	\label{Tab2}
\end{table}

In total, over $34,000$ numerical experiments were carried out according to the setup. Models were implemented using Julia and Gurobi version 7.5 on a Lenovo desktop machine with 8 GB RAM and Intel Core i5-65 CPU with 2.50GHz using Windows 10 OS 64-bit. In Gurobi, we have used a time limit of $9000s$ for each problem instance and optimality is achieved once the optimality gap is below $0.01\%$.

\subsection{Data}

We tested the discrete, polyhedral and stochastic models using network data instances taken from the online SNDlib library\footnote{See \url{http://sndlib.zib.de}}, see \cite{Orlowski2010}. The particular network data considered in this work is Germany-50 with $50$ nodes and $176$ directed arcs (we also included arcs in opposite directions). There are three levels of aggregation for real-world traffic measurement data available. These are one full day (in 5 minute intervals), one full month (in 1 day intervals) and one whole year (in 1 month intervals).

\begin{figure}[tbp]
	\begin{minipage} {.48\textwidth}
		\centering
		\includegraphics[width=\linewidth]{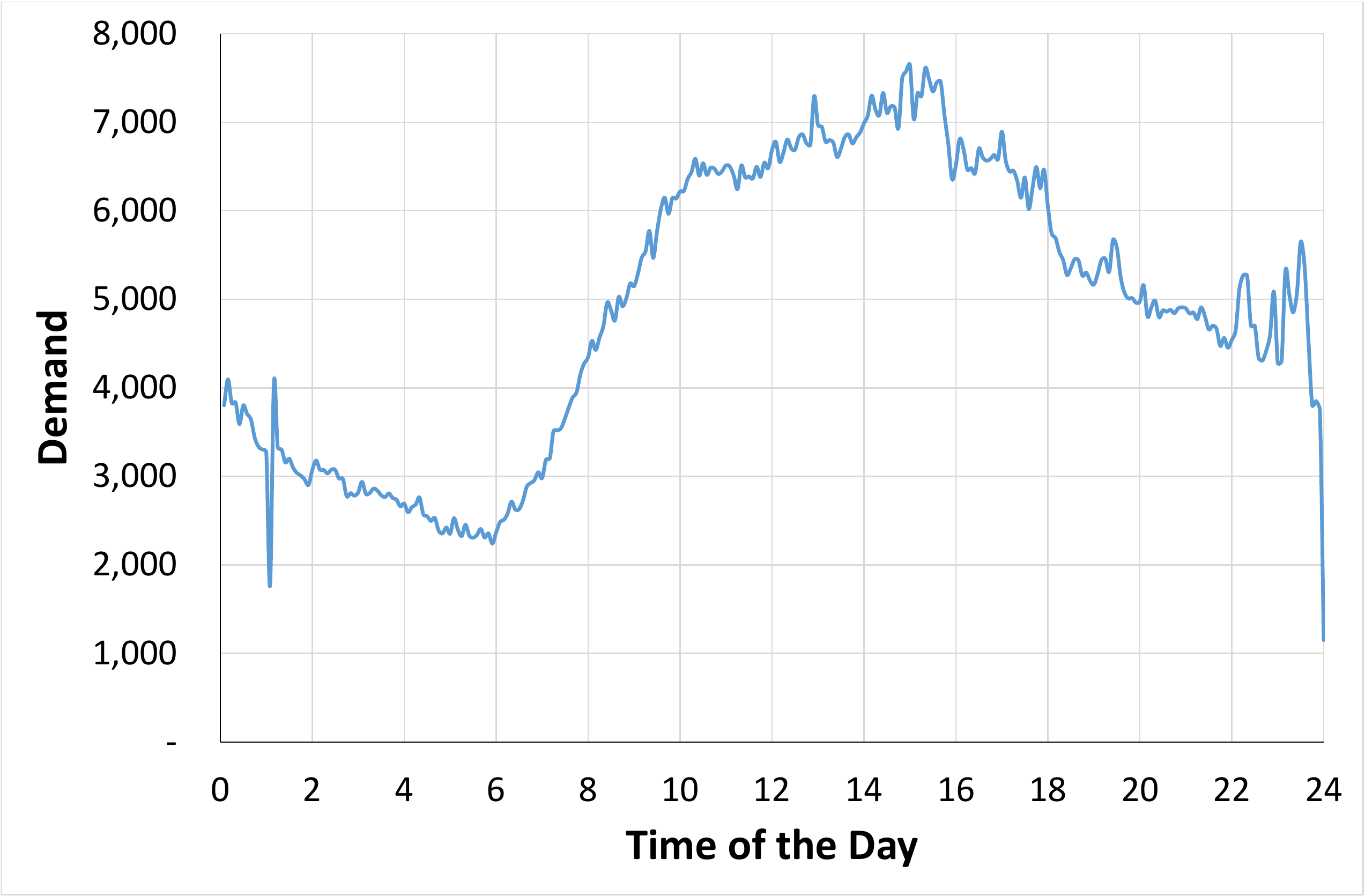}
		\caption{A full day demand profile.}\label{FigDP}
	\end{minipage}%
	\hfill
	\begin{minipage} {.48\textwidth}
		\centering
		\includegraphics[width=\linewidth]{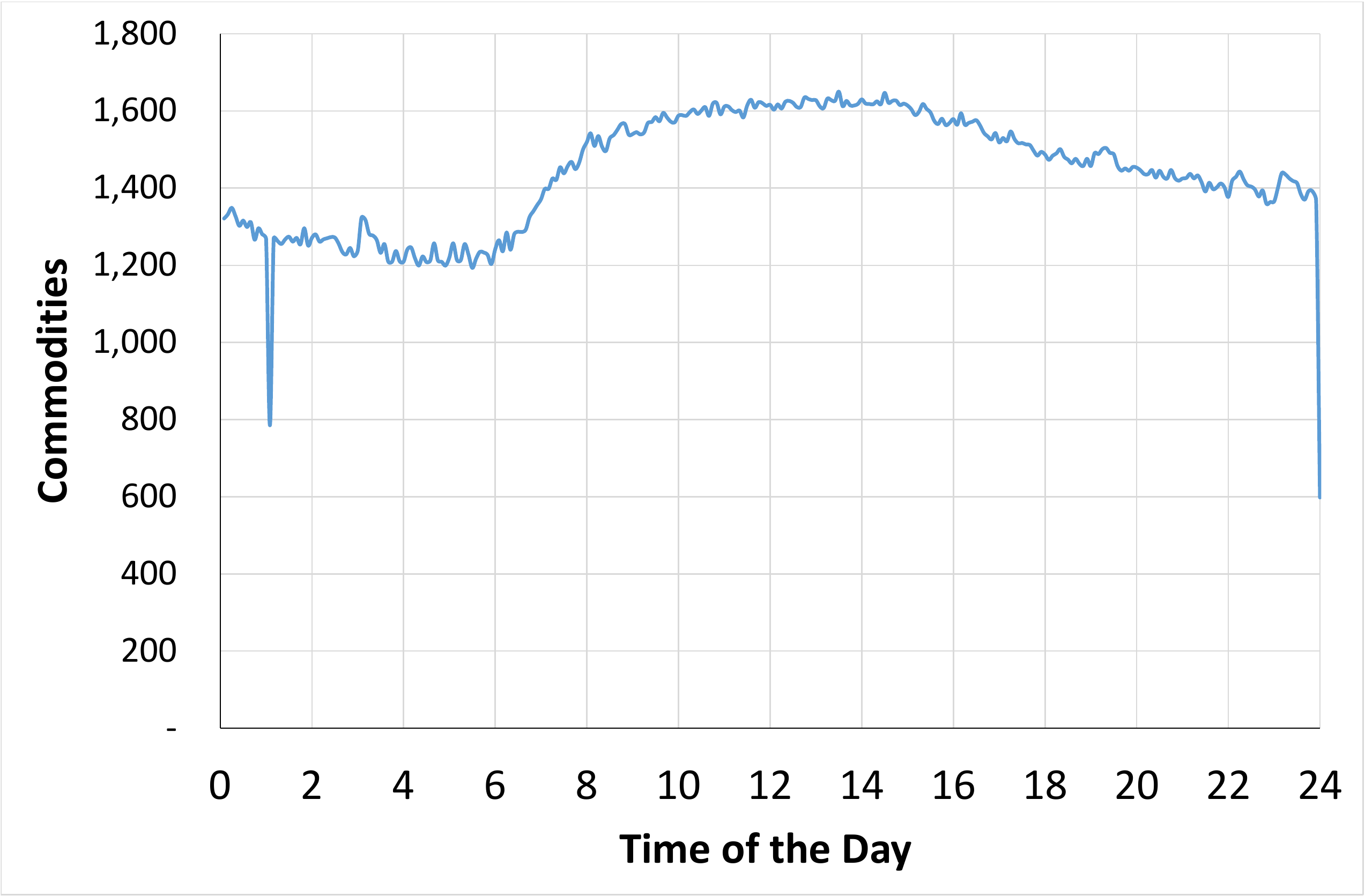}
		\caption{A full day commodities profile.}\label{FigCP}
	\end{minipage}	
\end{figure}

For our experiment, we focus on the full day dataset, consisting of $N=288$ scenarios.
The peak demand of $7,649.83$  was recorded at 3pm for the demand profile, see Figure~\ref{FigDP}.
We separate the scenarios into a training set consisting of 24 scenarios, which is generated by taking every 12th demand scenario (i.e., one scenario per hour), and the evaluation set consisting of the remaining 264 scenarios. We refer to the training set as MS-24.

Each scenario has a different number of commodities, see Figure~\ref{FigCP}.
Some of the demand values were found to be very small. While the $99^{th}$ percentile of all demand values is 0.415, some values are in the range of $10^{-6}$. To simplify the optimization problems, we sort the commodities in descending order of demand and then choose a fixed value of commodities for all demand scenarios that covers over 98\% of the original demand data, which is the case for 400 commodities. \autoref{Tab1} shows the different numbers of commodities against the percentage of original data captured in the streamlined data. This approach was implemented instead of allowing for varying commodities per demand scenario and allows us to consider all significant demands values while discarding very low ones, thus significantly reducing the average numbers of commodities per demand scenario.
\begin{table}[tbp]\footnotesize
	\centering
	%\tbl{Impact of Choice of k on the presorted original data.}
	\caption{Impact of choice of $K$ on the presorted original data.}
	\begin{tabular}{lcc}%ccc}
		\toprule
		\textbf{Options with MS-24} & \textbf{Commodity} &\textbf{\% of Original Data Captured} \\
		\toprule %\hline \hline
		All Demand & 300 & 97.48\% \\
		All Demand & 400 & 98.88\% \\
		All Demand & 450 & 99.25\% \\%$ %\widehat{d}$ \& $0.6 \widehat{d}$\\ %\hline
		All Demand & 500 & 99.50\% \\
% 		All Demand & 598 & 99.76\% \\
		All Demand & 900 & 99.88\% \\
% 		Demand Cut off of 0.4633 & 381 & 98.93\% \\
% 		Demand Cut off of 0.4134 & 393 & 99.04\% \\
		\bottomrule
	\end{tabular}
	%\caption{Experimental setup for generating 120 problem instances for each network.}
	\label{Tab1}
\end{table}

We observed that a model based on a polyhedron with 400 commodities computed from the training set demand matrix could not be solved in reasonable time, hence the polyhedron was generated for a reduced number of commodities  to allow for an optimal solution in a reasonable amount of time that will encourage its practical usage in the industry.
Instead, we work with $K=20$ that captures the top commodities in the training set. This reflects that the more complex the model for the uncertainty, the harder becomes the optimization model itself, and the less data we can use for building our sets. This trade-off is investigated in our experiments. Additionally, our polyhedrons were calcualted using the random constraint sampling method from Section~\ref{polcon}, as lower and upper bounds already gave an optimal solution to the optimization approach for constructing polyhedra.

\subsection{Computational Results}

We consider the performance of the capacity expansion solutions on the evaluation scenarios. We used four metrics on these 264 scenarios: The average, the maximum, the average of the worst 10\% (known as conditional-value-at-risk, or CVaR), and the standard deviation. Note that all these measures were calculated for scenarios that were not known to the models at the time of solution.

\begin{figure}[tbp]
	\begin{minipage}[t]{.48\textwidth}
		\centering
		\includegraphics[width=\textwidth]{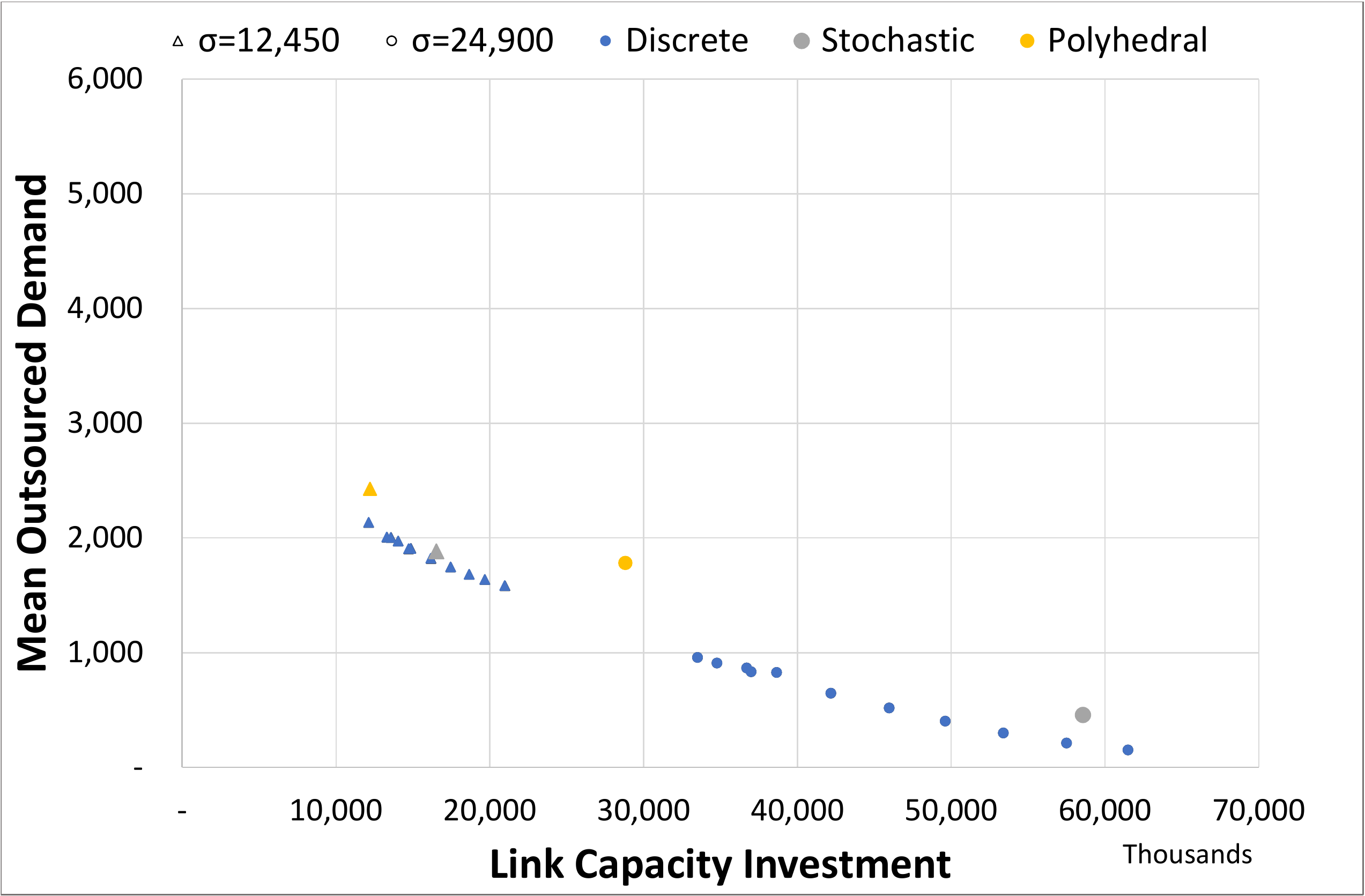}
		\caption{Mean outsourced demand. Discrete model uses varying values of $\lambda$.}\label{Fig1}
	\end{minipage}%
	\hfill
	\begin{minipage}[t]{.48\textwidth}
		%\centering
		\includegraphics[width=\textwidth]{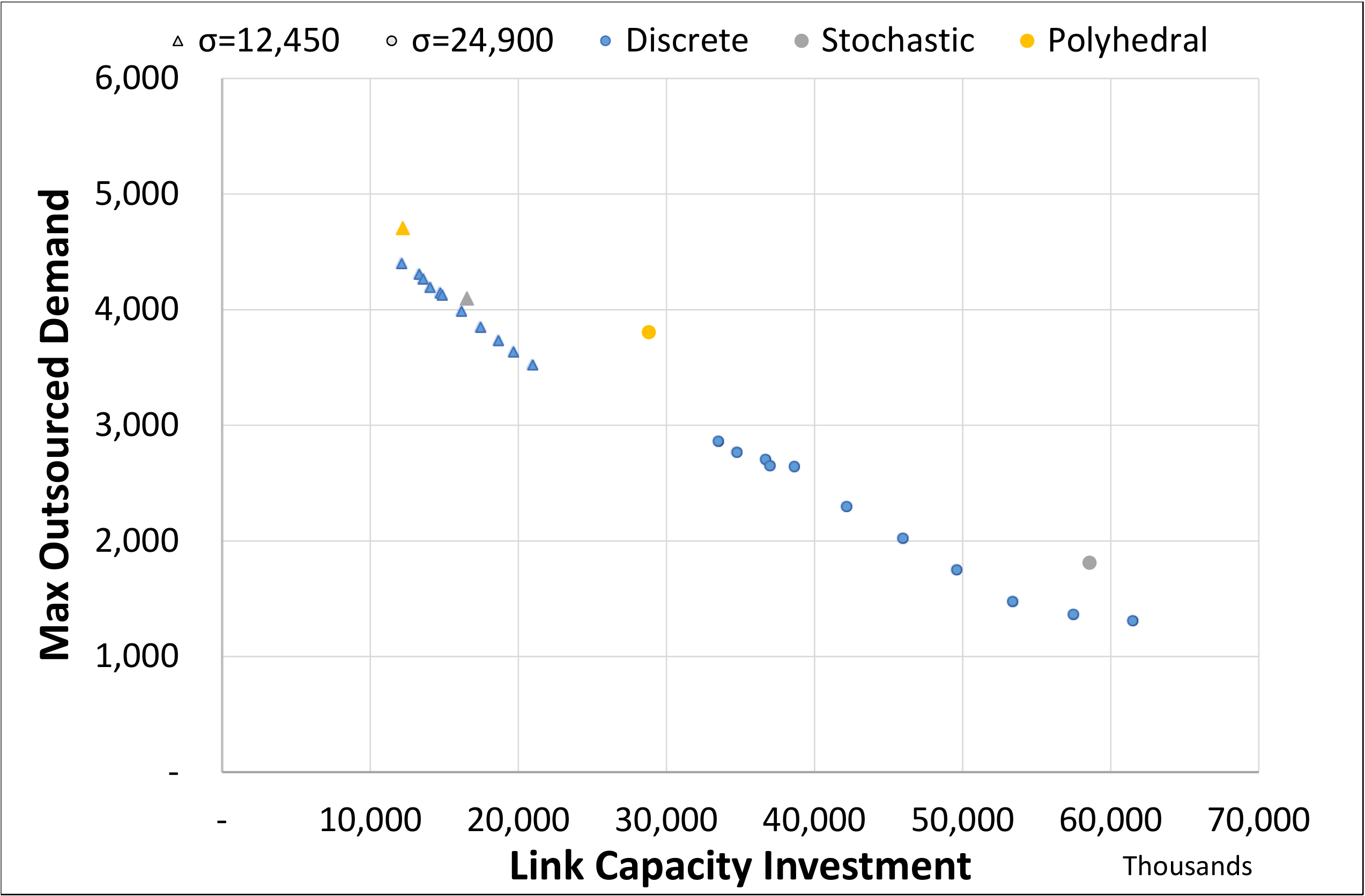}
		\caption{Maximum outsourced demand. Discrete model uses varying values of $\lambda$.}\label{Fig2}
	\end{minipage}	
\end{figure}

\begin{figure}[tbp]
	\begin{minipage}[t] {.48\textwidth}
		\centering
		\includegraphics[width=\linewidth]{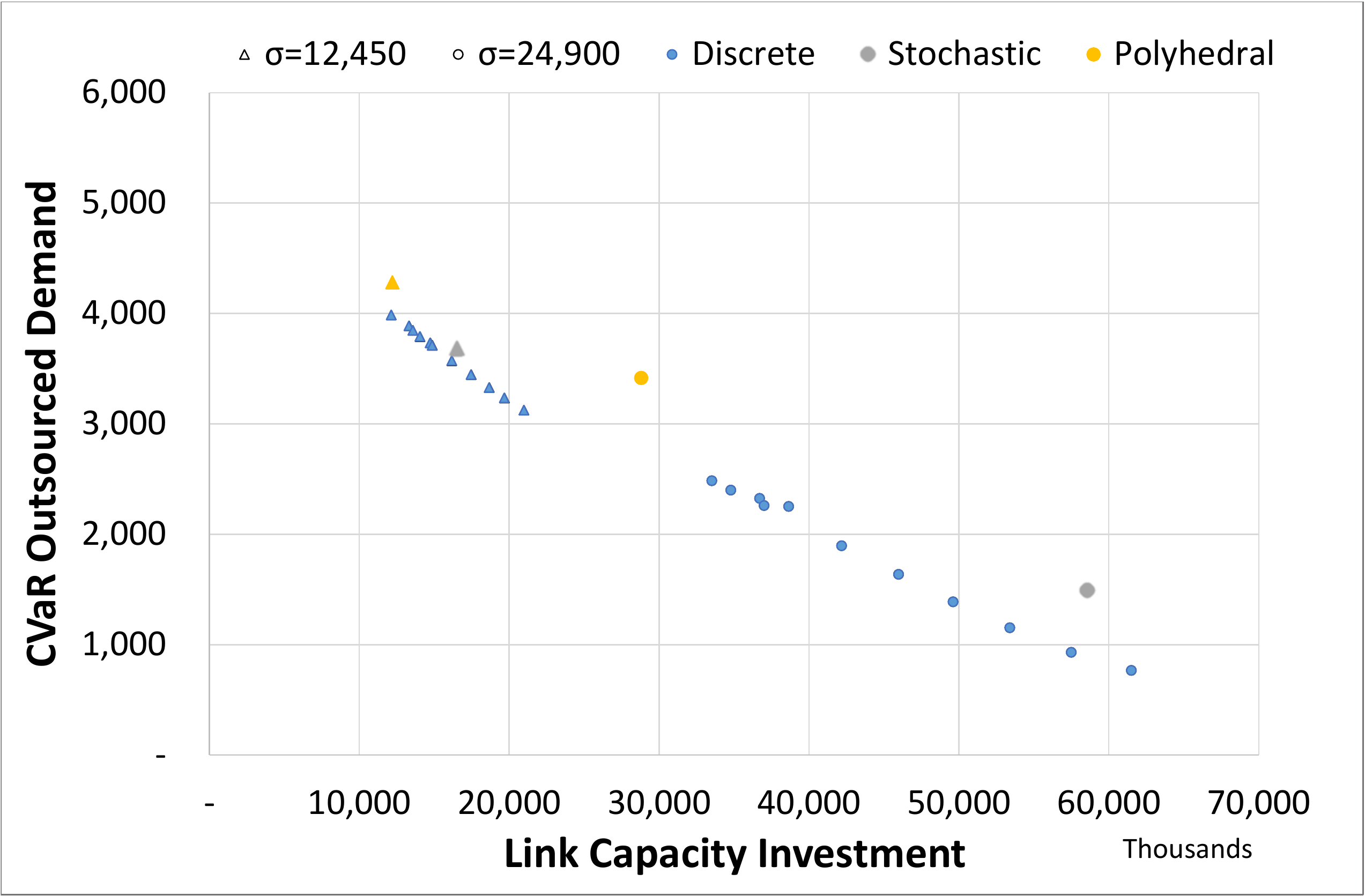}
		\caption{CVaR of outsourced demand. Discrete model uses varying values of $\lambda$.}\label{Fig3}
	\end{minipage}%
	\hfill
	\begin{minipage}[t] {.48\textwidth}
		%\centering
		\includegraphics[width=\linewidth]{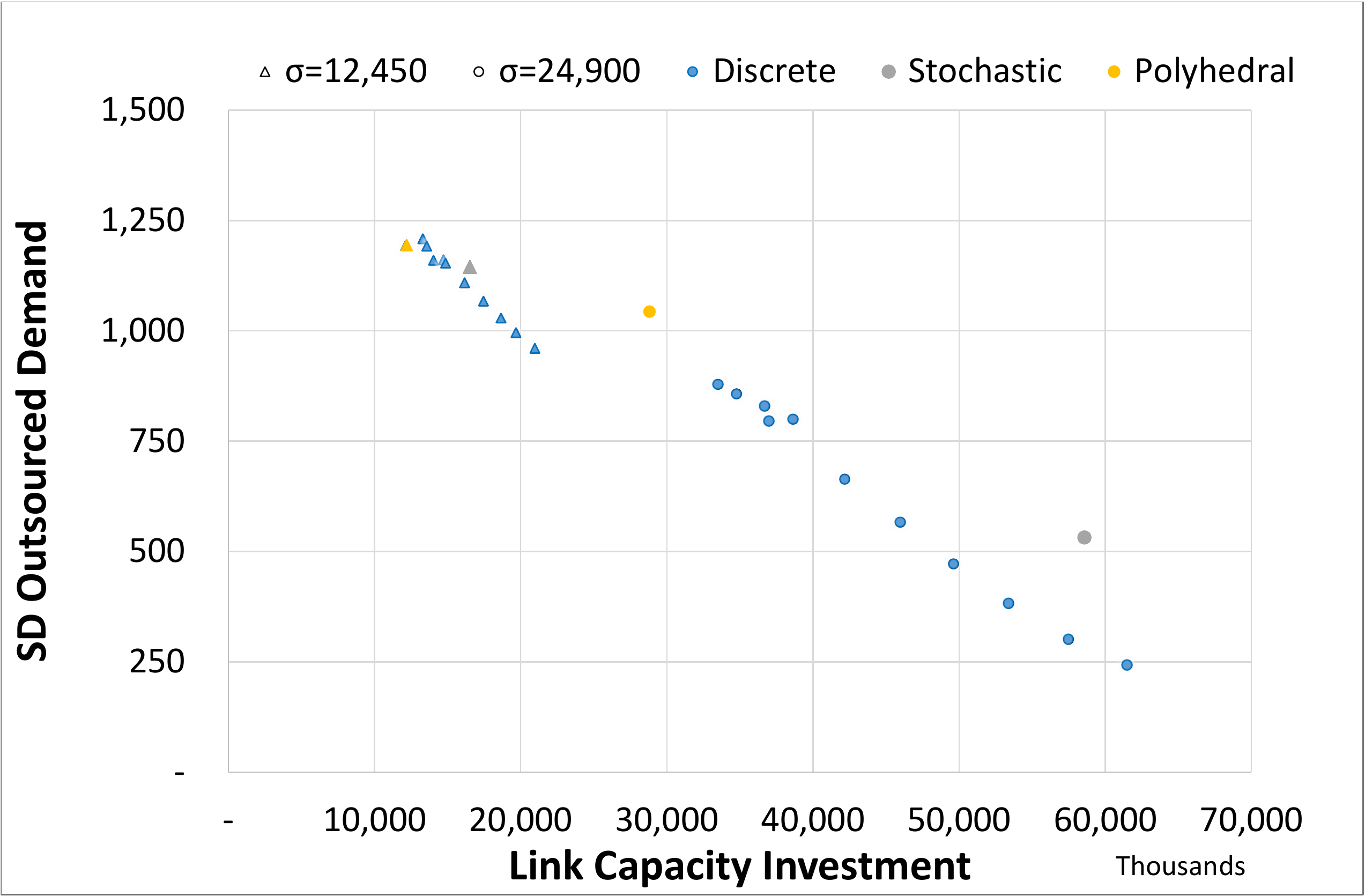}
		\caption{Standard deviation of outsourced Demand. Discrete model uses varying values of $\lambda$.}\label{Fig4}
	\end{minipage}	
\end{figure}

We first of all note that all polyhedral models $\text{Polyhedral}_1$ to $\text{Polyhedral}_3$ gave the same results, so we do not differentiate between them in the following. In \autoref{Fig1} to  \autoref{Fig4}, the four metrics are shown against the first-stage investment costs for two values of $\sigma$ (i.e., using $\text{Discrete}_1$). As expected, increasing $\sigma$ results in building more capacity in the network and hence reducing the amount of demand being outsourced. This is true for both robust the stochastic models. For the discrete uncertainties, network capacity built increases with increasing value of $\lambda$ from 0 (ignoring uncertainty) to 1 (using the real demands) for a fixed $\sigma$ value.

\begin{figure}[tbp]
	\begin{minipage}[t]{.48\textwidth}
		\centering
		\includegraphics[width=\linewidth]{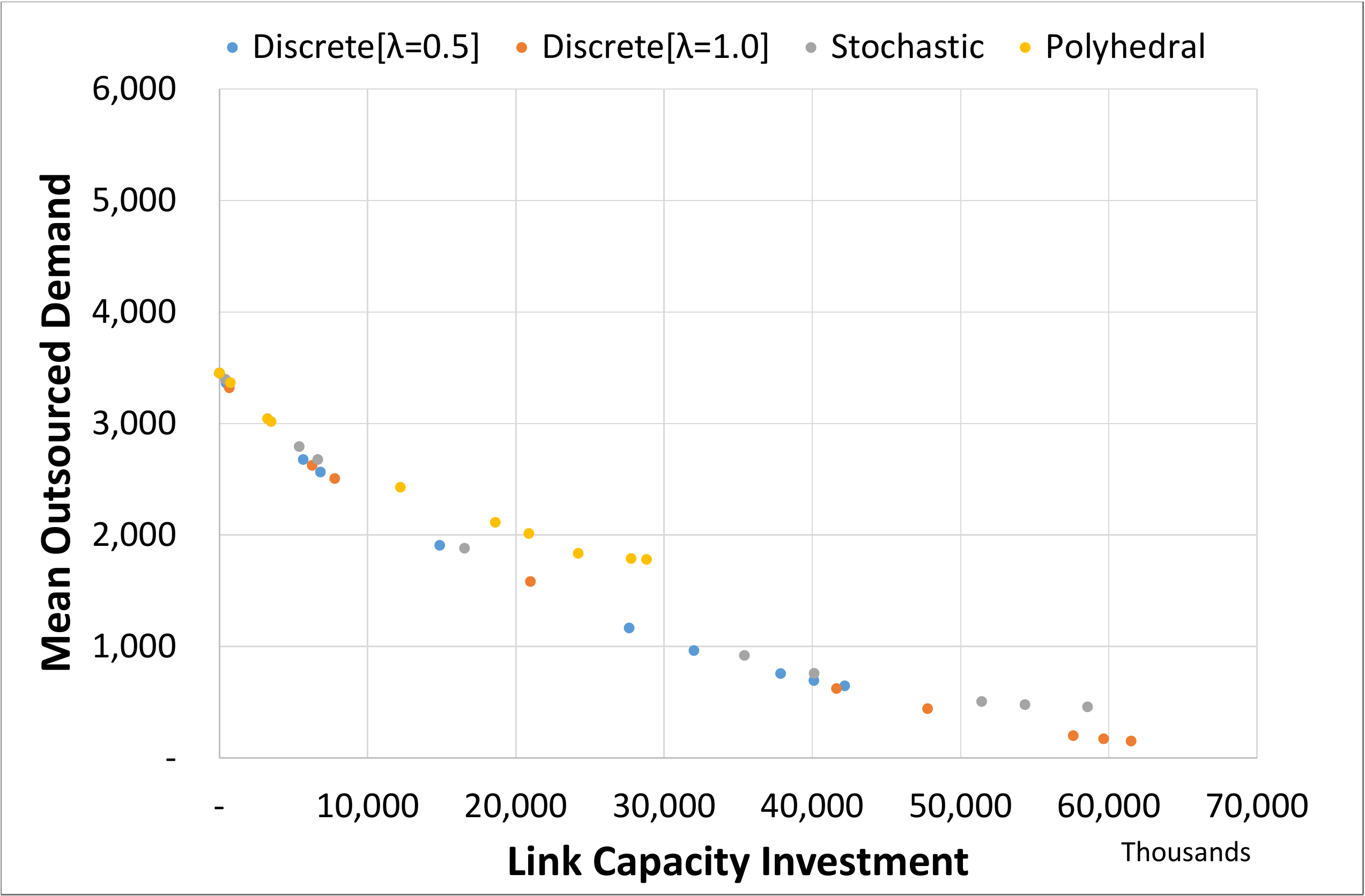}
		\caption{Mean outsourced demand. All models use varying values of $\sigma$.}\label{Fig5}
	\end{minipage}%
	\hfill
	\begin{minipage}[t]{.48\textwidth}
		%\centering
		\includegraphics[width=\linewidth]{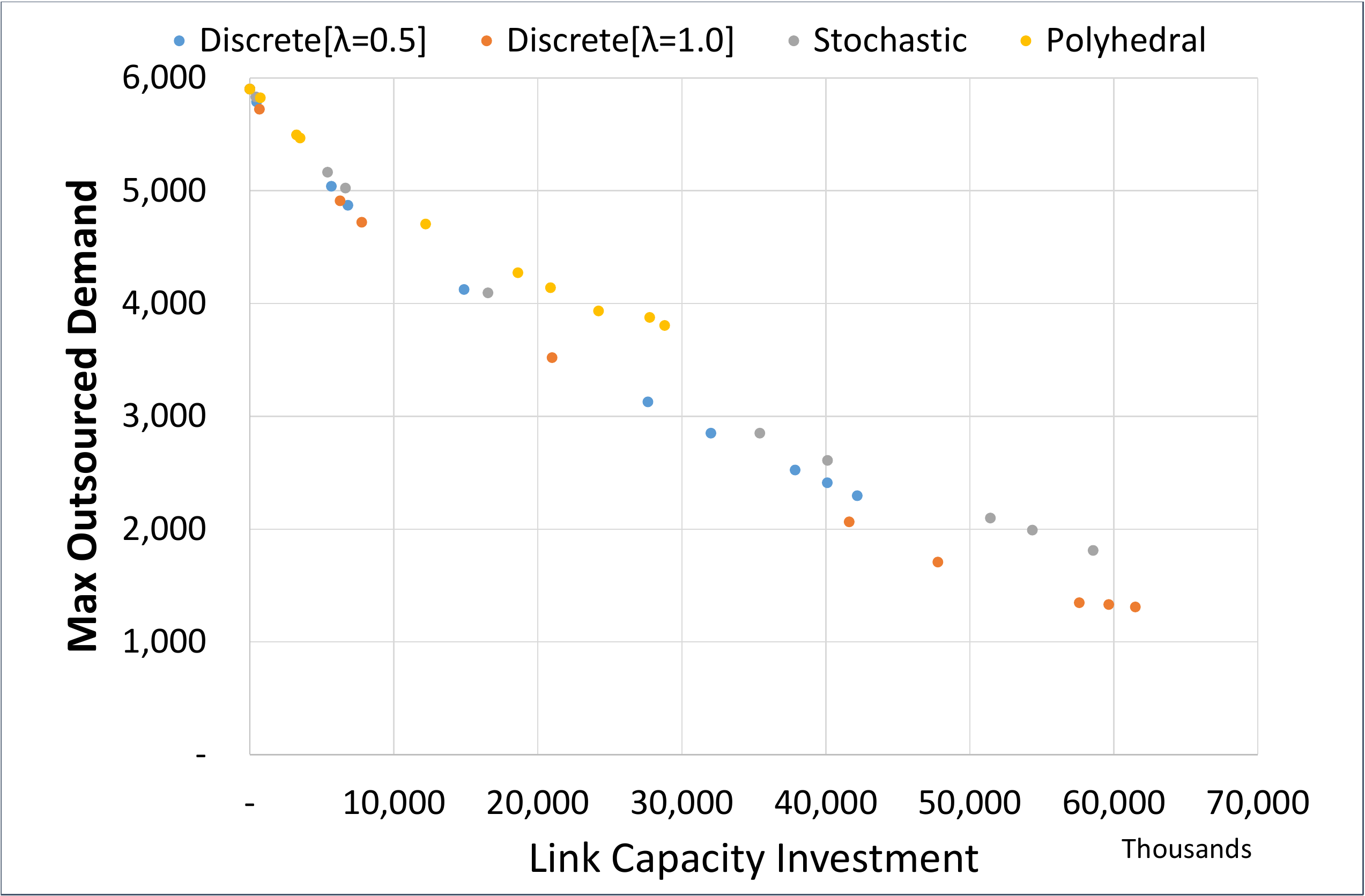}
		\caption{Maximum outsourced demand. All models use varying values of $\sigma$.}\label{Fig6}
	\end{minipage}	
\end{figure}

\begin{figure}[tbp]
	\begin{minipage}[t] {.48\textwidth}
		\centering
		\includegraphics[width=\linewidth]{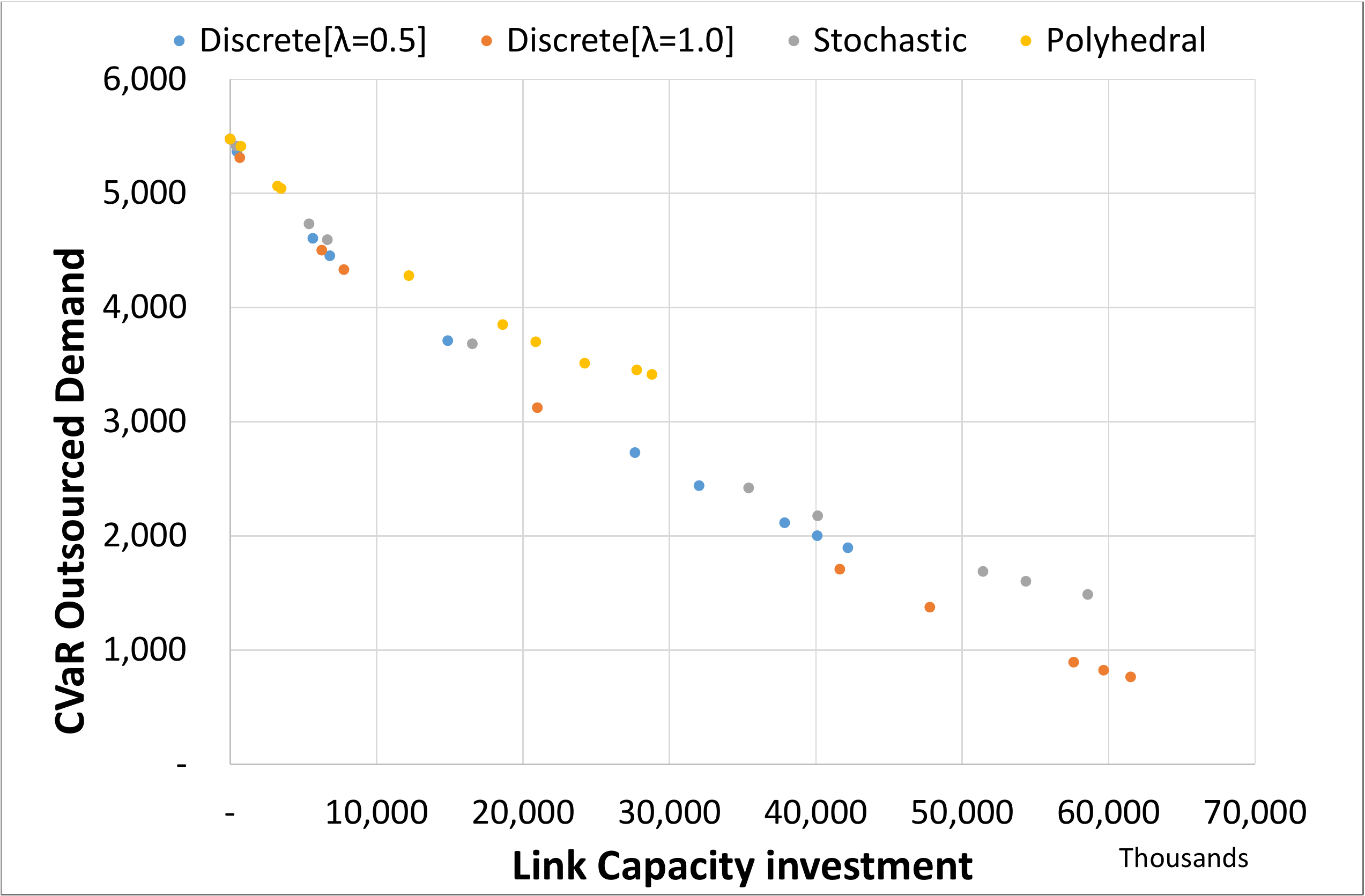}
		\caption{CVaR of outsourced demand. All models use varying values of $\sigma$.}\label{Fig7}
	\end{minipage}%
	\hfill
	\begin{minipage}[t] {.48\textwidth}
		%\centering
		\includegraphics[width=\linewidth]{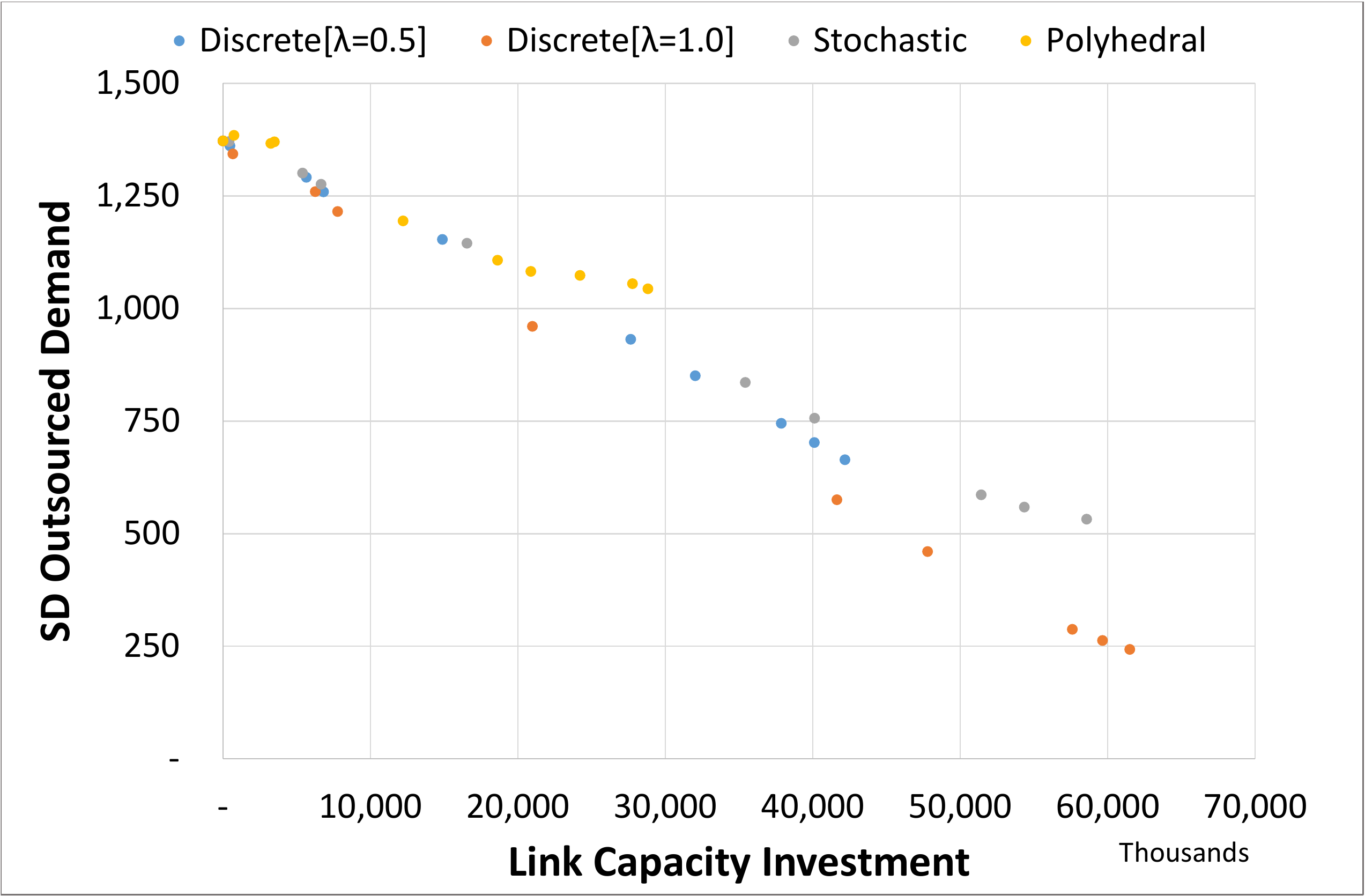}
		\caption{Standard deviation of outsourced demand. All models use varying values of $\sigma$.}\label{Fig8}
	\end{minipage}	
\end{figure}

\begin{figure}[tbp]
% 	\begin{minipage} {.48\textwidth}
		\centering
		\includegraphics[width=.48\linewidth]{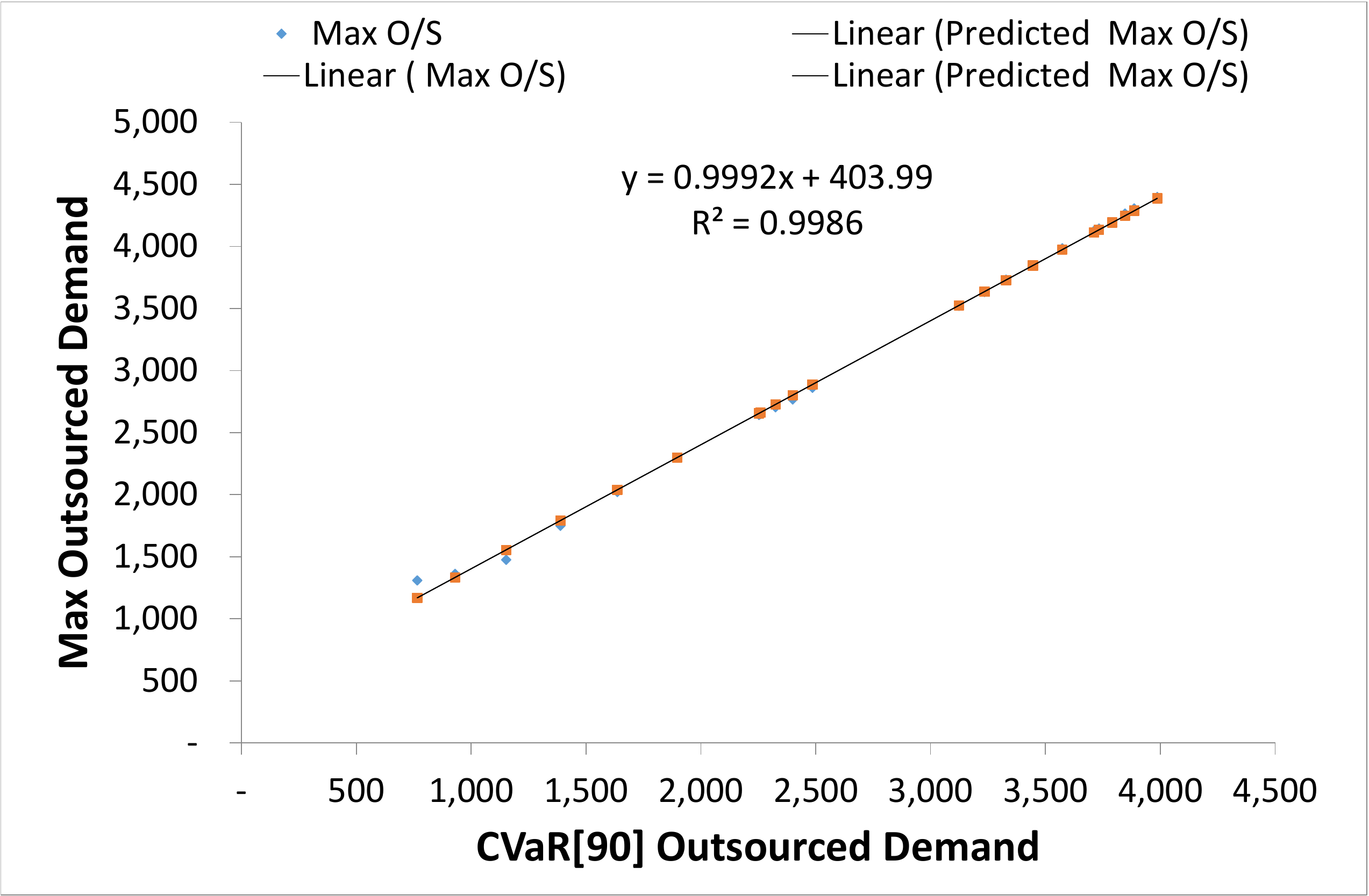}
		\caption{CVaR of outsourced demand and max outsourced demand correlation.}\label{Fig9}
% 	\end{minipage}%
% 	\begin{minipage} {.5\textwidth}
% 		%\centering
% 		\includegraphics[width=1.1 \linewidth]{"max OD si"}
% 		\caption{Maximum Outsourced Demand ($\sigma$).}\label{Fig10}
% 	\end{minipage}	
\end{figure}

In \autoref{Fig5} to \autoref{Fig8}, varying penalty values $\sigma$ were considered for the three models while fixing $\lambda$ for the discrete uncertainty model (using $\text{Discrete}_2$). The outsourced demand $\tau$ decreases with an increase in $\sigma$ value.
The implication of higher penalty is that overall risk is minimized deploying additional infrastructure in capacity for the network rather than outsourcing demand. In \autoref{Fig3} and \autoref{Fig7}, the CVaR was observed to decrease with increasing robustness of the models. \autoref{Fig2} seems to be providing almost the same information as the CVaR, and it turned out that the two metrics are highly correlated having a correlation coefficient of $0.9993$ with a gradient of approximately $1$ as shown in \autoref{Fig9}. Though the analysis done was for the discrete model, the same result is consistence with that from the other two models.

Ideally, a good solution is in the bottom left corner of these plots. We note that some of the points corresponding to polyhedral models are dominated, and so are the stochastic solutions.
The discrete model produces the best trade-off solutions between investment and outsourcing. For instance in \autoref{Fig7}, with the same link capacity investment of $\$40$ million, the stochastic model has a higher CVaR figure. The data point line for this discrete model with $\lambda=0.5$ is below that for the stochastic model and this can bee seen in \autoref{Fig5} to \autoref{Fig8}. Hence, the discrete model provides the best compromise between a too simple and a too complex approach for the data under consideration.

\section{Conclusion}
\label{sec:conclusion}

In the robust optimization literature, the shape of uncertainty is often an assumption made without any grounding in actually available data. This also holds for network expansion problems, where polyhedral models have been popular. In this paper, we considered the question whether such an approach leads to solutions which perform well on unseen data, i.e., what kind of uncertainty sets are most appropriate for our model.

We developed robust (using discrete and polyhedral uncertainty sets) and stochastic approaches to a multi-commodity network capacity expansion problem with the option of demand outsourcing. These models were implemented for a real-world network data taken from the SNDlib and their results were subsequently compared.

In the experimental setup, a number of penalty values for demand outsourcing were considered while also varying the robustness of the discrete model with
different sizes of the uncertainty set. Increasing the penalty results in additional capital expenditure for network capacity build as this reduces the
amount of demand outsourced as well as the conditional-value-at-risk (CVaR). However, of these three models, the robust model with discrete uncertainty set
produced the best trade-off solutions on all performance metrics. It was also observed that the discrete set seems easy to generate (as expected, since the
original data is already in this form), the model is simple and produces optimal result faster. Robust model with polyhedral uncertainty set, on the other
hand, is more complex and with more options to describe data, and it results in computationally more challenging problems. In our case, the extra effort
associated with polyhedral model may not be really worth it in the end. Surprisingly, the simple stochastic optimization model which we have used for
benchmarking was relatively competitive, and thus might be appropriate for use in more complex situations in which the uncertainty-based robust models are
computationally intractable.

\end{document}